\documentclass[12pt]{article}
\usepackage[dvips]{graphicx}
\DeclareGraphicsExtensions{ps,eps}
\makeatletter
\def\section{\@startsection{section}{1}{\z@}{-1.5ex plus -.5ex
minus -2.ex}{1ex plus .2ex}{\large\bf}}
\input{amssym.def}
\input{amssym}

\def\@thmcounterstep{}
\long\def\@makecaption#1#2{\vskip 10pt \setbox\@tempboxa\hbox{#1.#2}
\ifdim \wd\@tempboxa >\hsize
#1.#2\par 
\else99
\hbox to\hsize{\hfil\box\@tempboxa\hfil}
\fi}
\def\ps@headings{
\def\@oddhead{\footnotesize\rm\hfill\runninghead\hfill}
\def\@evenhead{\@oddhead}
\def\@oddfoot{\rm\hfill\thepage\hfill}\def\@evenfoot{\@oddfoot}}
\newtheorem{Theorem}{Theorem}[section]

\newtheorem{Proposition}[Theorem]{Proposition}

\makeatother
\title
{
HIDDEN CONVEXITY IN SOME NONLINEAR 
PDEs FROM GEOMETRY AND PHYSICS
}
\def\runninghead{\quad 
Hidden convexity in nonlinear PDEs
}
\author{
{\em Yann Brenier}\thanks{CNRS, Universit\'e de Nice (FR 2800 W. D\"oblin),
Institut Universitaire de France
}}
\date{} % no date wanted.
\begin{document}
\pagestyle{headings}
\flushbottom
\maketitle
\vspace{-10pt}%include according to taste.
%\tiny

\section{Introduction}

There is a prejudice among  some specialists of non linear partial differential
equations and differential geometry: convex analysis is an elegant theory
but too rigid to address some of the most interesting and challenging problems
in their field.
Convex analysis is mostly attached to elliptic
and parabolic equations  of variational origin,
for which a suitable convex potential can be exhibited and shown to be
minimized (either statically or dynamically).
The Dirichlet principle for linear elliptic equation is archetypal.
\\
Hyperbolic PDEs, for example, seem to be inaccessible to convex analysis,
since they are usually derived from variational principles that are definitely
not convex.  (However, convexity plays an important role in the so-called
entropy conditions.) Also, elliptic systems with variational
formulations (such as in elasticity theory)
often involve structural conditions quite far from convexity (such as
Hadamard's "rank one" conditions). (However, convexity can be 
often restored, for example through the concept of
polyconvexity \cite{Ba}, or by various kinds of "relaxation" methods
\cite{Yo,ABM}.)
The purpose of the present paper is to show few examples of nonlinear
PDEs (mostly with strong geometric features) for which there is a hidden
convex structure. This is not only a matter of curiosity. Once the
convex structure is unrevealed, robust existence and uniqueness results
can be unexpectedly obtained for very general data.
Of course, as usual,
regularity issues are left over as a hard post-process, but, at least,
existence and uniqueness  results are obtained in a large framework.
The paper will address:
\begin{enumerate}
\item 
{THE MONGE-AMPERE EQUATION}\\
(solving the Minkowski problem and strongly related to the so-called
optimal transport theory since the 1990's)
\item 
{THE EULER EQUATION}\\
(describing the motion of inviscid and incompressible fluids,
interpreted by Arnold as geodesic curves on infinite dimensional groups of
volume preserving diffeomorphisms)
\item 
{THE MULTIDIMENSIONAL HYPERBOLIC SCALAR CONSERVATION LAWS}
(a simplified model for multidimensional systems of
hyperbolic conservation laws)
\item 
{THE BORN-INFELD SYSTEM}\\
(a non-linear electromagnetic model introduced in 1934,
playing an important role in high energy Physics since the 1990's)
\end{enumerate}

Finally, let us mention that we borrowed the expression ``hidden convexity''
from a lecture by L.C. Evans about various models where the same
phenomena occur (such as growing sandpiles
\cite{AEW} and weak KAM theory).

\section{Monge-Amp\`ere equation and optimal transportation maps}

Given two positive functions $\alpha$ and $\beta$ of same 
integral over $\bf{R^d}$, we look for a convex solution $\Phi$ of the
the Monge-Amp\`ere equation:
\begin{equation}
\label{monge ampere}
\beta(\nabla\Phi(x))det(D^2\Phi(x))=\alpha(x),\;\;x\in R^d.
\end{equation}
This nonlinear PDE is usually related to the
Minkowski problem, which amounts to find hypersurfaces of prescribed
Gaussian curvature.

\subsection{A weak formulation}

Assuming $a$ $priori$ that $x\in R^d\rightarrow \nabla\Phi(x)$
is a diffeomorphism (with a jacobian matrix $D^2\Phi(x)$
everywhere symmetric positive), we immediately see, using the change of
variable $y=\nabla\Phi(x)$, that
(\ref{monge ampere}) is equivalent to
the following "weak formulation":
\begin{equation}
\label{weak monge ampere}
\int f(y)\beta(y)dy=\int f(\nabla\Phi(x)))\alpha(x)dx
\end{equation}
for all suitable test function $f$ on $R^d$.
In the words of measure theory, this just means that 
$\beta(y)dy$ as a Borel measure on $R^d$ is the image of
the measure $\alpha(x)dx$ by the map $x\rightarrow \nabla\Phi(x)$.
Notice that such a weak formulation has nothing to do with the usual
definition of weak solutions in the sense of distribution (that does not make
sense for a fully non-linear equation such as  (\ref{monge ampere})).
It is also weaker than the concept of ``viscosity solution'', as
discussed in \cite{Ca}.

\subsection{A convex variational principle}

\begin{Proposition}
\label{minimization}
Let us consider all smooth convex functions $\Psi$ on $R^d$ with
a smooth Legendre-Fenchel transform
\begin{equation}
\label{legendre}
\Psi^*(y)=\sup_{x\in R^d}x\cdot y -\Psi(x).
\end{equation}
Then, in this family, a solution $\Phi$ to the Monge-Amp\`ere 
equation (\ref{monge ampere}) is a minimizer of the convex functional
\begin{equation}
\label{functional}
J[\Psi]=
\int \Psi(x)\alpha(x)dx+\int \Psi^*(y)\beta(y)dy.
\end{equation}

\end{Proposition}

\subsubsection*{Proof}
For any suitable convex function $\Psi$, we have:
$$
J[\Psi]=\int \Psi(x)\alpha(x)dx+\int \Psi^*(y)\beta(y)dy
=\int (\Psi(x)+\Psi^*(\nabla\Phi(x)))\alpha(x)dx
$$
(since $\nabla\Phi$ transports $\alpha$ toward $\beta$)
$$
{\ge \int x\cdot\nabla\Phi(x)\alpha(x)dx}
$$
(by definition of the Legendre transform (\ref{legendre}))
$$
= \int (\Phi(x)+\Phi^*(\nabla\Phi(x)))\alpha(x)dx
$$
(indeed, in the definition of ${\Phi^*(y)=\sup \;\;x\cdot y -\Phi(x)}$,
the supremum is achieved whenever ${y=\nabla\Phi(x)}$,
which implies ${\Phi^*(\nabla\Phi(x))=x\cdot\nabla\Phi(x)-\Phi(x)}$)
$$
=\int \Phi(x)\alpha(x)dx+\int \Phi^*(y)\beta(y)dy=J[\Phi],
$$
which shows that, indeed, $\Phi$ is a minimizer for (\ref{functional}))

\subsection{Existence and uniqueness for the weak Monge-Amp\`ere problem}

Based on the previous observation,  using the tools of convex analysis,
one can solve the Monge-Amp\`ere problem in its weak formulation,
for a quite large class of data, with both existence and uniqueness of
a solution:

\begin{Theorem}
\label{OT}
Whenever ${\alpha}$ and ${\beta}$ 
are nonnegative Lebesgue integrable functions on $R^d$,
with same integral,
and bounded second order moments,
$$
{\int |x|^2\alpha(x)dx<+\infty,\;\;\;\int |y|^2\beta(y)dy<+\infty,}
$$
there is a unique $L^2$  map $T$ with convex potential 
$ T=\nabla\Phi$
that solves the Monge-Amp\`ere problem in its weak formulation
(\ref{weak monge ampere}), for all continuous function $f$
such that: $|f(x)|\le 1+ |x|^2$.
 \\
This map is called the optimal transport map between
$\alpha(x)dx$ and $\beta(y)dy$.
\end{Theorem}

By map with convex potential, we exactly means a Borel map $T$
with the following property:
there is a lsc convex
function $\Phi$ defined on $R^d$, valued in $]-\infty,+\infty]$,
such that, for $\alpha(x)dx$ almost everywhere $x\in R^d$,
$\Phi$ is differentiable at $x$ and $\nabla\Phi(x)=T(x)$.

\subsubsection*{Comments}

The usual proof  \cite{Br3,SK,RR,Br5} is based on the duality method introduced
by  Kantorovich to solve the so-called Monge-Kantorovich problem,
based on the concept of joint measure (or coupling measure) \cite{Ka}. 
However a direct proof is possible, as observed by Gangbo \cite{Ga}.
This theorem 
can be seen as the starting point of the
so-called "optimal transport theory" which has turned out to be a very
important and active field of research in the recent years, with a lot of interactions
between calculus of variations, convex analysis, differential geometry,
PDEs, functional analysis and probability theory and several
applications outside of mathematics (see \cite{Vi} for a review).
A typical (and striking) application to the isoperimetric inequality 
is given in the appendix of the present paper.

\section{ The Euler equations}

\subsection{ Geometric definition of the Euler equations}

The Euler equations were introduced in 1755 \cite{Eu}
to describe the motion of inviscid  fluids.
In the special case of  an incompressible fluid 
moving inside a bounded convex domain $D$ in $R^d$,
a natural configuration space is the set
$SDiff(D)$ of all orientation and volume preserving diffeomorphisms of $D$.
Then, a solution of the Euler equations can be defined as a curve
$t\rightarrow g_t$ along $SDiff(D)$ subject to:
\begin{equation}
\label{euler}
\frac{d^2 g_t}{dt^2}\circ g_t^{-1}+\nabla p_t=0,
\end{equation}
where
$p_t$ is a time dependent scalar field defined on $D$ (called the 'pressure field').
\\
As shown by Arnold \cite{AK,EM}, these equations have a very simple geometric
interpretation.  Indeed,  $g_t$ is just a geodesic curve (with constant speed)
along $SDiff(D)$, with respect to the $L^2$
metric inherited from the Euclidean space $L^2(D,R^d)$, and $-\nabla p_t$ is the acceleration
term, taking into account the curvature of $SDiff(D)$.
From this interpretation in terms of geodesics, we immediately deduce a  variational principle for the Euler equations.  However, this principle cannot be convex due to the non convexity of
the configuration space. (Observe that $SDiff(D)$ is contained in a sphere of the space
$L^2(D,R^d)$ and cannot be convex, except in the trivial case $d=1$ where it reduces to
the identity map.)

\subsection{A concave maximization principle for the pressure}

Surprisingly enough, the pressure field obeys (at least on short time intervals)
a concave maximization principle. More precisely,

\begin{Theorem}
\label{pressure maximization}
Let $(g_t,p_t)$ a smooth solution to the Euler equations (\ref{euler})
on a time interval $[t_0,t_1]$ small enough so that
\begin{equation}
\label{smallness}
(t_1-t_0)^2 D^2 p_t(x)\le \pi^2, \;\;\;\forall x\in D
\end{equation}
(in the sense of symmetric matrices).  Then $p_t$ is a maximizer
of the $CONCAVE$ functional
\begin{equation}
\label{functional2}
q\;\;\Rightarrow \;\;\int_{t_0}^{t_1}\int_D q_t(x)dtdx
+\int_D J_q[g_{t_0}(x),g_{t_1}(x)]dx,
\end{equation}
among all $t$ dependent scalar field $q_t$ defined on $D$.
Here
\begin{equation}
\label{subfunctional}
J_q[x,y]=\inf \int_{t_0}^{t_1} 
(-q_t(z(t))+\frac{|z'(t)|^2}{2})dt,
\end{equation}
where the infimum is taken over all curves ${t\rightarrow z(t)\in D}$
such that ${z(t_0)=x\in D}$, ${z(t_1)=y\in D}$, is defined for all pair of points $(x,y)$ in $D$.
\end{Theorem}

\subsubsection*{Proof}
The proof  is very elementary and
does not essentially differ from the one we used for the Monge-Amp\`ere equation in the
previous section
(which is somewhat surprising since the Euler equations and the MA equation look quite
different). The main difference is the smallness condition we need on the size of the
time interval.
Let us consider a  time dependent scalar field $q_t$ defined on $D$.
By definition of ${J_q}$:
$$
\int_D J_q[g_{t_0}(x),g_{t_1}(x)]dx
\le \int_{t_0}^{t_1}\int_D 
(\frac{1}{2}|\frac{dg_t}{dt}|^2-q_t(g_t(x)))dtdx.
$$
Using a standard variational argument, we see that,  under the smallness condition (\ref{smallness}), the Euler equation (\ref{euler}) 
asserts that, for all $x\in D$
$$
 J_p[g_{t_0}(x),g_{t_1}(x)]
=\int _{t_0}^{t_1}
(\frac{1}{2}|\frac{dg_t}{dt}|^2-p_t(g_t(x)))dt.
$$
Integrating in $x\in D$, we get:
$$
\int_D J_p[g_{t_0}(x),g_{t_1}(x)]dx
= \int_{t_0}^{t_1}\int_D 
(\frac{1}{2}|\frac{dg_t}{dt}|^2-p_t(g_t(x)))dtdx.
$$
Since $g_t\in SDiff(D)$ is volume preserving, we have:
$$
\int_D (q_t(x)-q_t(g_t(x))dx=\int_D (p_t(x)-p_t(g_t(x))dx=0.
$$
Finally,
$$
\int_{t_0}^{t_1}\int_D q_t(x)dtdx
+\int_D J_q[g_{t_0}(x),g_{t_1}(x)]dx
$$
$$
\le 
\int_{t_0}^{t_1}\int_D p_t(x)dtdx
+\int_D J_p[g_{t_0}(x),g_{t_1}(x)]dx
$$
which shows that, indeed, $(p_t)$ is a maximizer.

\subsection{
Global convex analysis of the Euler equations
}

The maximization principle is the starting point for a global analysis of the Euler equations.
Of course,  there is no attempt here to solve
the Cauchy problem in the large for $d\ge 3$, which is one of the most outstanding
problems in nonlinear PDEs theory. (This would more or less
amount to prove the geodesic completeness of $SDiff(D)$.)
We rather address the existence of minimizing geodesics between 
arbitrarily given points of the configuration space $SDiff(D)$. This problem may have
no classical solution, as shown by Shnirelman \cite{Sh1}.
Combining various contributions by Shnirelman, Ambrosio-Figalli and the author
\cite{Br4,Sh2,Br6,AF}, we get the global existence and uniqueness result:

\begin{Theorem}
\label{euler global}
Let  $g_0$ and $g_1$ 
be given volume preserving Borel maps of $D$ (not necessarily
diffeomorphisms) and $t_0<t_1$. Then
\\
1) There is a ${unique}$  $t$ dependent pressure field $p_t$, with zero mean on
$D$, that solves 
(in a suitable weak sense) the maximization problem stated in Theorem  \ref{pressure maximization}
\\
2) There is a sequence ${{g^n_t}}$ 
valued in  SDiff(D) such that
$${\frac{d^2 g^n_t}{dt^2}\circ (g_t^n)^{-1}+\nabla p_t\rightarrow 0,}$$
in the sense of distributions and
${g^n_0\rightarrow g_0,\;\;\;g^n_1\rightarrow g_1}$ in ${{L^2}}$.
\\
3) Any sequence of approximate minimizing geodesics $(g^n_t)$ (in a suitable sense)
betwween $g_0$ and $g_1$ has the previous behaviour.
\\
4) The pressure field is well defined in  the space $L^2(]t_0,t_1[,BV_{loc}(D))$.

\end{Theorem}

Of course, these results are not as straightforward as Theorem  \ref{pressure maximization} and
requires a lot of technicalities (generalized flows, etc...).
However, they still rely on
convex analysis which is very surprising in this infinite 
dimensional differential
geometric setting. Notice that the uniqueness result is also surprising. 
Indeed, between two given points, minimizing
geodesics are not necessarily unique (as can be easily checked). However the
corresponding acceleration field $-\nabla p_t$ is unique! It is unlikely that such a property
could be proven using classical differential geometric tools. It is  probably an output of  
the hidden convex structure.
Let us finally notice that the improved regularity obtained by Ambrosio
and Figalli \cite{AF} (they show that $p$ belongs to
$L^2(]t_0,t_1[,BV_{loc}(D))$ instead of $\nabla p$ a locally bounded measure, 
as previously obtained in \cite{Br6}) is just sufficient to give a full
meaning to the maximization problem. (A different formulation, involving
a kind of Kantorovich duality is used in \cite{Br4,Br6} and requires less
regularity.)

\section{ Convex formulation of  multidimensional scalar conservation laws}
\subsection{Hyperbolic systems of conservation laws}

The general form of multidimensional nonlinear conservation laws is:
$$
\partial_t u_t+\sum_{i=1}^d \partial_i (F_i(u_t))=0,
$$
where $u_t(x)\in V\subset R^m$ is a time dependent vector-valued field defined on a $d-$ dimensional
manifold (say the flat torus  $T^d=R^d/Z^d$ for simplicity)
and each $F_i:V\subset R^m\rightarrow R^m$ is a given nonlinear function.
This general form includes systems of paramount importance in Mechanics and Physics,
such as the gas dynamics and the Magnetohydrodynamics equations, for example.
A simple necessary (and nearly sufficient) condition for the Cauchy problem to be 
well-posed for short times 
is the hyperbolicity condition which requires, for all $\xi\in R^d$ and all $v\in V$
the $m\times m$ real matric 
$$
\sum_{i=1}^d \xi_i F_i'(v)
$$
to be diagonalizable with real eigenvalues. For many systems of physical origin, with a
variational origin, there is
an additional conservation law:
\begin{equation}
\label{entropy system}
\partial_t (U(u_t))+\sum_{i=1}^d \partial_i (G_i(u_t))=0,
\end{equation}
where $U$ and $G_i$ are scalar functions (depending on $F$). (This usually follows
from Noether's invariance theorem.) Whenever, $U$
is a strictly convex function, the system automatically gets hyperbolic.
For most  hyperbolic systems, solutions are expected to become discontinuous
in finite time, even for smooth initial conditions. There is no  theory available to
solve the initial value problem in the large
(see \cite{Da} for a modern review), except in two extreme situations.
First, for a single space variable ($d=1$)
and small initial conditions (in total variation), global existence and
uniqueness of "entropy solutions" have been established through the celebrated results of
J. Glimm (existence) and A. Bressan and collaborators (well posedness) 
\cite{Gl,BB}. (Note that some special systems can also be treated 
with the help of compensated compactness methods \cite{Ta}, without
restriction on the size of the initial conditions.)
Next, in the multidimensional case, global existence and uniqueness of 'entropy solutions'  have been
obtained by Kruzhkov \cite{Kr} in the case of a single (scalar) conservation law ($m=1$).

\begin{Theorem}
\label{kruzhkov}
(Kruzhkov)
\\
Assume $F$ to be Lipschitz continuous.
Then, for all $u_0\in L^1(T^d)$, there is a unique $(u_t)$,
in the space $C^0(R_+,L^1(T^d))$ with initial value $u_0$, such that:
\\
\begin{equation}
\label{scalar}
\partial_t u_t+\nabla \cdot (F(u_t))=0,
\end{equation}
is satisfied in the distributional sense and,
\\
for all Lipschitz convex function $U$ defined on  $R$, the "entropy" inequality
\begin{equation}
\label{entropy}
\partial_t (U(u_t))+\nabla \cdot (Z(u_t))\le 0,
\end{equation}
holds true in the distributional sense, where
$$
Z(v)=\int_0^v F'(w)U'(w)dw.
$$
In addition, for all pair of such "entropy" solutions $(u,\tilde u)$,
\begin{equation}
\label{contraction}
\int_{T^d}  |u_t(x)-\tilde u_t(x)|\;dx\;\le\;\int_{T^d}  |u_s(x)-\tilde u_s(x)|\;dx,
\;\;\;\forall t\ge s\ge 0.
\end{equation}
\end{Theorem}

This result is often quoted as a typical example of maximal monotone operator theory in $L^1$. (For the concept of maximal monotone operator, we refer
to \cite{Brz,ABM}.)
The use of the non hilbertian space $L^1$ is crucial.
Indeed (except in the trivial linear case $F(v)=v$),
the entropy solutions do not depend on their initial values
in a  Lipschitz continuous way in any space $L^p$ except fot $p=1$. This is due to the fact that,
even for a smooth initial condition, the corresponding entropy solution $u_t$ may
become discontinuous for some $t>0$ and, therefore, cannot belong to any Sobolev space 
$W^{1,p}(T^d)$ for $p>1$.

\subsection{A purely convex formulation of multidimensional scalar conservation laws}

Clearly, convexity is already involved in Kruzhkov's formulation 
(\ref{scalar},\ref{entropy})
of scalar conservation laws, through the concept of "entropy inequality". 
However, a deeper, hidden, convex structure
can be exhibited, as observed recently by the author \cite{Br10}.
As a matter of fact,  the Kruzhkov entropy solutions can be fully recovered just by
solving a rather straightforward convex sudifferential inequality in the Hilbert space
$L^2$. For notational  simplicity, we limit ourself to the case when the initial
condition $u_0$ is valued in the unit interval.

\begin{Theorem}
(YB, 2006, \cite{Br10})
\\
Assume $u_0(x)$ to be valued in $[0,1]$, for $x\in T^d$.
Let $Y_0(x,a)$ be any bounded function of $x\in T^d$ and $a\in [0,1]$,
non decreasing in $a$,  such that
\begin{equation}
\label{expansion0}
u_0(x)=\int_0^1 1\{Y_0(x,a)<0\}da,
\end{equation}
for instance: $Y_0(a)=a-u_0(x)$.
Then, the unique Kruzhkov solution to (\ref{scalar}) is given by
\begin{equation}
\label{expansion}
u_t(x)=\int_0^1 1\{Y_t(x,a)<0\}da,
\end{equation}
where $Y_t$ solves the convex
subdifferential inequality in $L^2(T^d\times [0,1])$:
\begin{equation}
\label{master}
0\; \epsilon \;\partial_t  Y_t+F'(a)\cdot\nabla_x Y_t+\partial \eta[Y_t],
\end{equation}
where $\eta[Y]=0$ if $\partial_a Y\ge 0$, 
and $\eta[Y]=+\infty$ otherwise.
\end{Theorem}

Observe that  $Y\rightarrow F'(a)\cdot\nabla_x Y+\partial \eta[Y]$ defines a maximal monotone
operator  and generates a semi-group of contractions in  $L^2(T^d\times [0,1])$  \cite{Brz}.

\subsubsection*{Sketch of proof}

Multidimensional scalar conservation laws enjoy a comparison principle (this is why they
are so simple with respect to general systems of conservation laws). In other words,
if a family of initial conditions $u_0(x,y)$ is non decreasing with respect to a real parameter $y$,
the corresponding Kruzhkov solutions $u_t(x,y)$ will satisfy the same property. This key
observation enables us to use a kind of level set method, in the spirit of Sethian and
Osher \cite{OS,OF}, and even more closely, in the spirit of the paper by Tsai, Giga and Osher \cite{TGO}.
Assume, for a while, that $u_t(x,y)$ is $a$ $priori$ smooth and strictly increasing with respect to $y$.
Thus, we can write
$$
u_t(x,Y_t(x,a))=a,\;\;\;Y_t(x,u_t(x,y))=y
$$
where $Y_t(x,a)$ is smooth and strictly increasing in $a\in [0,1]$.
Then, a straightforward calculation shows that $Y$ must solve the simple linear equation
\begin{equation}
\label{linear}
\partial_t  Y_t+F'(a)\cdot\nabla_x Y_t=0
\end{equation}
(which has $Y_t(x,a)=Y_0(x-tF'(a),a)$ as exact solution).
Unfortunately, this linear equation is not able to preserve the monotonicity condition
$\partial_a Y\ge 0$ in the large.  Subdifferential inequality (\ref{master}) is, therefore,
a natural substitute for it. The remarkable fact is that this rather straighforward 
modification exactly matches  the Kruzhkov entropy inequalities.
More precisely, as $Y$ solves (\ref{master}) ,  then
$$
u_t(x,y)=\int_0^1 1\{Y_t(x,a)<y\}da
$$
can be shown to be the right entropy solutions with initial conditions $u_0(x,y)$.
For more details, we refer to \cite{Br10}.

\subsubsection*{Remark}

Our approach is reminiscent of both 
the "kinetic method" and the "level set" method.
The kinetic approach amounts to linearize the scalar conservation laws as (\ref{linear}) 
by adding an extra variable (here $a$).  This idea (that has obvious roots in the
kinetic theory of Maxwell and Boltzmann)
was independently introduced for scalar conservation laws 
by Giga-Miyakawa and the author  \cite{Br1,Br2,GM}. Using this approach,
Lions, Perthame and Tadmor \cite{LPT}  later introduced the so-called
kinetic formulation of scalar conservation laws and, using the averaging lemma
of Golse, Perthame and Sentis  \cite{GPS}, established the remarkable result
that multidimensional scalar conservation laws
enjoy a regularizing effect when they are genuinely nonlinear. 
(In other words, due to shock waves, entropy solutions automatically get
a fractional amount of differentiability!).
On the other side, the level set method by Osher 
and Sethian \cite{OS,OF} describes functions according to their level sets (here $Y(t,x,a)=y$).
This is a very general and powerful approach to all kinds of numerical
and analytic issues in pure and applied mathematics.  An application of the level set
method to scalar conservation laws  was made by Tsai, Giga and Osher  in \cite{TGO}
and more or less amounts to introduce a viscous (parabolic) approximation of subdifferential
inequality (\ref{master}).
Finally, let us mention that some very special systems of conservation laws
can be treated in a similar way \cite{Br8}.

%%%%%%%%%%%%%%%%%%%%%%%%%%%%%%%%%%%%%%%%%%%%%%%%%%%%%% SLIDE19

\section{ The Born-Infeld system}

Using convential notations of classical electromagnetism, the Born-Infeld system reads:
$$
{{{
{\partial_t B+\nabla\times (\frac{B\times(D\times B)+D}
{\sqrt{1+D^2+B^2+(D\times B)^2}})=0,\;\;\;
\nabla\cdot B=0,}
}}}
$$
$$
{{{
{\partial_t D+\nabla\times (\frac{D\times(D\times B)-B}
{\sqrt{1+D^2+B^2+(D\times B)^2}})=0,\;\;\;
\nabla\cdot D=0,}
}}}
$$
This system is a nonlinear correction to the Maxwell equations,
which can describe strings and branes in high energy Physics 
\cite{Bo,BI,BDLL,Gi}.
Concerning the initial value problem,
global smooth solutions have been proven to exist for small localized
initial conditions by Chae and Huh \cite{CH}  (using Klainerman's null forms and following
a related work by Lindblad \cite{Li} ).
The additional conservation law
$$
{{{
{\partial_t h+\nabla\cdot Q=0,}
}}}
$$
where
$$
{{{
{h=\sqrt{1+D^2+B^2+(D\times B)^2},\;\;
Q=D\times B.}
}}}
$$
provides an 'entropy function' ${h}$ which is a convex function of the
unknown ${(D,B)}$ only in a neighborhood of $(0,0)$.
However, $h$ is clearly a convex function of $B$, $D$, and $B\times D$.
Thus, there is a hope to restore convexity by considering $B\times D$
as an independent variable, which will be done subsequently by
``augmenting'' the Born-Infeld system.

%%%%%%%%%%%%%%%%%%%%%%%%%%%%%%%%%%%%%%%%%%%%%%%%%%%%%% SLIDE6
%%%%%%%%%%%%%%%%%%%%%%%%%%%%%%%%%%%%%%%%%%%%%%%%%%%%%% SLIDE3

\subsection{ The augmented  Born-Infeld (ABI) system}
Using Noether's invariance theorem, we get from the BI system
$4$ additional ('momentum-energy' ) conservation laws:
\begin{equation}
\label{ABI}
\partial_tQ+ 
\nabla \cdot (\frac{Q\otimes Q
-B\otimes B-D\otimes D}{h})=\nabla(\frac{1}{h}),\;\;\;\;\;
\partial_t h+\nabla\cdot Q=0.
\end{equation}
We call augmented Born-Infeld system (ABI) the $10\times 10$ 
system of equations made of 
the $6$ original BI evolution equations
\begin{equation}
\label{BI}
\partial_t B+\nabla\times (\frac{B\times Q+D}{h})=
\partial_t D+\nabla\times (\frac{D\times Q-B}{h})=0,
\end{equation}
with the differential constraints
\begin{equation}
\label{differential}
\nabla\cdot B=0,\;\;\;\nabla\cdot D=0,
\end{equation}
together with the 4 additional conservation laws (\ref{ABI}) 
but WITHOUT the algebraic constraints
\begin{equation}
\label{manifold}
h=\sqrt{1+D^2+B^2+(D\times B)^2},\;\;
Q=D\times B.
\end{equation}
These algebraic constraints define a $6$ manifold in the space
$(h,Q,D,B)\in R^{10}$ that we call "BI manifold". 
We have the following consistency result:

\begin{Proposition}
(Y.B., 2004 \cite{Br7})
\\
Smooth solutions of the ABI system (\ref{ABI},\ref{BI},\ref{differential}) 
preserve the BI manifold (\ref{manifold}). Therefore, any smooth solution of
the original BI system can be seen as a smooth solution to the ABI system 
(\ref{ABI},\ref{BI},\ref{differential}) 
with an initial condition valued on the BI manifold.
\end{Proposition}

%%%%%%%%%%%%%%%%%%%%%%%%%%%%%%%%%%%%%%%%%%%%%%%%%%%%%% SLIDE6

\subsection{First appearance of convexity in the ABI system}

Surprisingly enough, the $10\times 10$ augmented  ABI system 
(\ref{ABI},\ref{BI},\ref{differential}) 
admits  an extra conservation law:
$$
\partial_t U+\nabla\cdot Z=0,
$$
where
$$
U(h,Q,D,B)=\frac{1+D^2+B^2+Q^2}{h}
$$
is convex (and $Z$ is a rational function of $h,Q,D,B$). This 
leads to the {GLOBAL} hyperbolicity of the system.
\\
Notice that the ABI system looks like Magnetohydrodynamics equations and enjoys 
${classical}$ Galilean invariance:
$$
(t,x)\rightarrow (t,x+u\;t),\;\;\;(h,Q,D,B)\rightarrow (h,Q-hu,D,B),
$$
for any constant speed ${u\in R^3}$!
\\
For a large class of nonlinear Maxwell equations, a similar extension
can be down (with 9 equations instead of 10) as in \cite{Se1}.
It should be mentioned that a similar method was introduced earlier
in the framework of nonlinear elastodynamics with polyconvex
energy (see \cite{Da}).
\subsection{Second appearance of convexity in the ABI system}

The $10\times 10$ ABI (augmented Born-Infeld) system
is $linearly$ $degenerate$ \cite{Da}
and stable under weak-* convergence:
weak limits of uniformly bounded sequences  in $L^\infty$ of smooth solutions
depending on one space variable only are still solutions.
(This can be proven by using the Murat-Tartar 'div-curl' lemma.)
Thus, we may conjecture that the convex-hull of the BI manifold is
a natural configuration space for the (extended) BI theory.
(As a matter of fact, the differential constraints 
$\nabla\cdot D=\nabla\cdot B=0$
must be carefully taken into account, as pointed out  to us by Felix Otto.)
The convex hull has full dimension. More precisely,
as shown by D. Serre \cite{Se2},  the convexified BI manifold is  just defined by 
the following inequality:
\begin{equation}
\label{hull}
h\ge \sqrt{1+D^2+B^2+Q^2+2|D\times B-Q|}.
\end{equation}
\\
Observe that, on this convexified BI manifold (\ref{hull}):
\\
1) The electromagnetic field ${(D,B)}$ and the 'density and momentum'
fields ${(h,Q)}$ can be chosen $independently$ of each other, 
as long as they
satisfy the required  $inequality$ (\ref{hull}). Thus,
in some sense, the ABI system describes a coupling between field and matter,
original Born-Infeld model is purely electromagnetic.
\\
2) 'Matter' may exist without electromagnetic field: ${B=D=0}$, 
which leads to
the Chaplygin gas (a possible model for 'dark energy' or 'vacuum energy')
$$
\partial_tQ+ 
\nabla \cdot (\frac{Q\otimes Q}{h})=\nabla(\frac{1}{h}),
\;\;\;\;\;
\partial_th+\nabla\cdot Q=0,
$$
\\
3) 'Moderate' Galilean transforms are allowed
$$
(t,x)\rightarrow (t,x+U\;t),\;\;\;(h,Q,D,B)\rightarrow (h,Q-hU,D,B),
$$
which is impossible on the original BI manifold (consistently with special
relativity) but becomes possible under weak completion (see the related
discussion on 'subrelativistic' conditions in \cite{Br9}).

\section{Appendix: A proof of the isoperimetric inequality using an optimal transport map}

In this appendix, we describe a typical and striking 
application of optimal transport map methods.
Let  $\Omega$ be a smooth bounded open set and $B_1$ the unit ball in $\bf{R}^d$.
The isoperimetric inequality reads (with obvious notations):
$$
{|\Omega|^{1-1/d}|B_1|^{1/d}\le \frac{1}{d}|\partial\Omega|}.
$$
Let ${\nabla\Phi}$ the optimal transportation map between 
$$
\alpha(x)=\frac{1}{|\Omega|}1\{x\in \Omega\},\;\;\;
\beta(y)=\frac{1}{|B_1|}1\{x\in B_1\}.
$$
In such a situation, according to  Caffarelli's regularity result \cite{Ca},  $\nabla\Phi$ is a
diffeomorphism between $\Omega$ and $B_1$ (up to their boundaries)
with $C^2$ internal regularity
(which is not a trivial fact) and
$$
det(D^2\Phi(x))=\frac{|B_1|}{|\Omega|}\;,\;\;\;x\in \Omega.
$$
holds true in the classical sense.
Then the proof (adaptated from Gromov) of the isoperimetric inequality
is straightforward and sharp. Indeed,
since  ${\nabla\Phi}$ maps $\Omega$ to the unit ball, we have:
$$
|\partial\Omega|=\int_{\partial\Omega}d\sigma(x)
\ge \int_{\partial\Omega} \nabla\Phi(x)\cdot n(x)d\sigma(x)
$$
(denoting by $d\sigma$ and $n(x)$ respectively the Hausdorff measure and
the unit normal along the boundary of $\Omega$)
$$
=\int_\Omega \Delta \Phi(x)dx
$$
(using Green's formula)
$$
{
\ge d\int_\Omega (det(D^2\Phi(x))^{1/d}dx}
$$
(using that ${(detA)^{1/d}\le 1/d\;\;Trace(A)}$
for any nonnegative symmetric matrix ${A}$)
$${=d|\Omega|^{1-1/d}|B_1|^{1/d}}\;\;$$ 
since
$
{\;\;det(D^2\Phi(x))=\frac{|B_1|}{|\Omega|}\;,\;\;\;x\in \Omega.}
$
\\
So, the isoperimetric inequality
$$
{|\Omega|^{1-1/d}|B_1|^{1/d}\le \frac{1}{d}|\partial\Omega|}
$$
follows, with equality  $only$ when ${\Omega}$ is a ball,
as  can be easily checked by tracing back the previous inequalities.
Notice that Gromov's original proof does not require the map $T$ to be optimal
(it is enough that its jacobian matrix has positive eigenvalues).
However,  the optimal map plays a crucial role for 
various refinements of the isoperimetric inequality
(in particular its quantitative versions by Figalli-Maggi-Pratelli \cite{FMP}, for example).

\end{document}